\documentclass{elsart3-1}

\usepackage{amssymb}
\usepackage{amsmath}
\usepackage[english,francais]{babel}

\newtheorem{theorem}{Theorem}[section]

\newtheorem{e-proposition}[theorem]{Proposition}

\newtheorem{e-definition}[theorem]{Definition\rm}

\newtheorem{theoreme}{Th\'eor\`eme}[section]
\newtheorem{lemme}[theoreme]{Lemme}

\renewcommand{\theequation}{\arabic{equation}}
\def\HH{{\cal H}}

\def\rr{\mathbb R}

\def\nn{\mathbb N}

\def\HH{\mathcal H}
\def\AA{\mathcal A}

\providecommand{\abs}[1]{\lvert#1\rvert}%absolute value
\providecommand{\norm}[1]{\lVert#1\rVert}%norm
\def\og{\leavevmode\raise.3ex\hbox{$\scriptscriptstyle\langle\!\langle$~}}
\def\fg{\leavevmode\raise.3ex\hbox{~$\!\scriptscriptstyle\,\rangle\!\rangle$}}

\journal{the Acad\'emie des sciences}
\begin{document}
% place in the next line the header (rubrique) chosen for your article,
% if you know it (you can also have 2, format : Header1/Header2
\centerline{}
\begin{frontmatter}


% \author{Name\thanksref{label2}}
\selectlanguage{english}
\title{Stabilisation faible interne locale de syst\`eme \'elastique de Bresse}

% The [label1] can be suppressed if there is only one address for all authors

\selectlanguage{english}
\author[a,b]{Nahla Noun},
\ead{nahlanoun@hotmail.com}
\author[a]{Ali Wehbe},
\ead{ali.wehbe@ul.edu.lb}
\medskip

\address[a]{Universit\'e Libanaise, Facult\'e des Sciences 1 \& Hadath,  Beyrouth, Liban}
\address[b]{Universit\'e Montpellier 2, ACSIOM, Place Eug\`ene Bataillon, 34095 Montpellier, France}

\medskip

\begin{abstract}
\selectlanguage{english}
% Text of abstract in English
{\bf Weakly locally internal stabilization of elastic Bresse system}\\
In \cite{AlabauBresse}, Alabau et al. studied the exponential and polynomial stability of the Bresse system with
one globally distributed dissipation law.
In this note, our goal is to extend the results from \cite{AlabauBresse},
by taking into consideration the important case when the dissipation law is locally distributed and to improve the polynomial energy decay rate.
We then study the energy decay rate of the Bresse system with one locally internal distributed dissipation
law  acting on the equation about the shear angle displacement. Under the equal speed wave propagation condition, we show that the system is exponentially stable. On the contrary, we establish a new polynomial energy decay rate.

\vskip 0.5\baselineskip

\selectlanguage{francais}
% Text of abstract in French
\noindent{\bf R\'esum\'e} \vskip 0.5\baselineskip
Dans \cite{AlabauBresse}, Alabau et al. ont \'etudi\'e la stabilisation exponentielle et polynomiale de syst\`eme de Bresse sous l'action d'une seule loi de dissipation
globalement distribu\'ee. Dans cette note, notre but est d'\'etendre les r\'esultats de \cite{AlabauBresse}, pour prendre en consid\'eration le cas important o\`u la loi de dissipation est localement distribu\'ee et pour am\'eliorer le taux de d\'ecroissance polynomial de l'\'energie.
Nous \'etudions alors, le taux de d\'ecroissance de l'\'energie
du syst\`eme de Bresse sous l'action d'une seule loi de dissipation interne localement distribu\'ee et agissant sur l'\'equation de rotation angulaire.
Sous la condition d'\'egalit\'e des vitesses de propagation, nous montrons que le syst\`eme est exponentiellement stable.
Dans le cas contraire, nous \'etablissons un nouveau taux de d\'ecroissance polynomial de l'\'energie.

\end{abstract}
\end{frontmatter}
\selectlanguage{english}
\vspace{0.1cm}

\noindent{\bf I- Abridged English Version}

\section{Introduction}\label{intro}
\setcounter{equation}{0} In this note, we study the locally internal stabilization of an elastic Bresse system with one damping acting on the equation about the shear angle displacement. We consider the system:
\begin{eqnarray}
\rho_1\varphi_{tt}-\kappa(\varphi_x+\psi+l\omega)_x-\kappa_0 l(\omega_x-l\varphi)&=&0\hspace{1 cm} {\rm in}\ \ (0,L)\times (0,\infty),\label{e11}\\
\rho_2\psi_{tt}-b\psi_{xx}+\kappa(\varphi_x+\psi+l\omega)+a(x)\psi_t&=&0\hspace{1 cm} {\rm in}\ \ (0,L)\times (0,\infty),\label{e12}\\
\rho_1\omega_{tt}-\kappa_0(\omega_x-l\varphi)_x+\kappa l(\varphi_x+\psi+l\omega)&=&0\hspace{1cm} {\rm in}\ \ (0,L)\times (0,\infty),\label{e13}
\end{eqnarray}
with the following boundary conditions
\begin{eqnarray}
\varphi(t,x)=\psi_x(t,x)=\omega_x(t,x)&=&0 \hspace{1 cm} {\rm for}\ \ x=0,L,\label{e15}\\
\varphi(t,x)=\psi(t,x)=\omega(t,x)&=&0 \hspace{1 cm} {\rm for}\ \  x=0,L,\label{e16}
\end{eqnarray}
where $\varphi$, $\psi$, $\omega$ are the vertical, shear angle and longitudinal displacements; $a\in L^{\infty}([0,L];\rr^+)$ and we assume that there exist $\alpha$, $\beta$ such that $0\leq\alpha<\beta\leq L$ and $a \geq a_0 >0$ on $]\alpha,\beta[$. Here $\rho_1=\rho A$, $\rho_2=\rho I$, $\kappa_0=EA$, $\kappa=\kappa^\prime GA$, $b=EI$ and $l=R^{-1}$
are positive constants for the elastic material properties. To be more precise, $\rho$ for density, $E$ for the modulus of elasticity,
$G$ for the shear modulus, $\kappa^\prime$ for the shear factor, $A$ for the cross-sectional area, $I$ for the second moment of area of cross-section, and $R$ for the radius of the curvature (for more details see Lagnese et al. \cite{Lagnese1}).

The Bresse system consists of three coupled wave equations, there are number of publications concerning the stabilization of this system \cite{LiuRao1}, \cite{WehbeBresse}, \cite{RiveraBresse} and \cite{AlabauBresse}. In particular, in \cite{AlabauBresse}, Alabau et al. consider the Bresse system (\ref{e11})-(\ref{e16}) with one globally dissipation law acting on the equation about the shear angle displacement i.e the function $a$ is constant
on $]0,L[$. Under the equal speed wave propagation condition, $\kappa=\kappa_0$ and $\frac{\rho_1}{\rho_2}=\frac{\kappa}{b}$,
they established an exponential energy decay rate for usual initial data. On the contrary, when $\kappa\not=\kappa_0$ or $\frac{\rho_1}{\rho_2}\not=\frac{\kappa}{b}$, they, first, showed that the Bresse system loses the exponential stability.
Next, under the condition $\kappa=\kappa_0$ and $\frac{\rho_1}{\rho_2}\not =\frac{\kappa}{b}$, they established a polynomial energy decay
rate of type $\frac{1}{t^{2/3}}$ for smooth initial data. Otherwise, they showed that the energy of the regular solution
decays polynomially to zero with rates $\frac{1}{t^{1/3}}$. Then, the stabilization of
the Bresse system with one locally internal distributed dissipation law and the optimality of
the polynomial decay rate seems still open problem. Therefore, our goal is to extend the results from \cite{AlabauBresse},
by taking into consideration the important case when $a$ is a function in $L^\infty([0,L];\rr_+)$ strictly positive in an open interval $]\alpha,\beta[\subset ]0,L[$ (the cases $\alpha=0$ or $\beta=L$ are not excluded) and to improve the polynomial energy decay rate.

In this note, we consider the Bresse system damped by one locally internal distributed dissipation
law with Dirichlet-Neumann-Neumann or Dirichlet-Dirichlet-Dirichlet boundary conditions type. Under the equal speed wave propagation condition,
$\kappa=\kappa_0$ and $\frac{\rho_1}{\rho_2}=\frac{\kappa}{b}$, we establish the same exponential energy decay rate for usual initial data.
On the contrary, when $\kappa\not=\kappa_0$ or $\frac{\rho_1}{\rho_2}\not=\frac{\kappa}{b}$, we first prove the non-exponential decay rate for the Bresse system with Dirichlet-Neumann-Neumann condition type. Therefore, under the condition $\kappa =\kappa_0$ and $\frac{\rho_1}{\rho_2}\not=\frac{\kappa}{b}$,
we establish a new polynomial energy decay rate of type $\frac{1}{t}$ for the smooth solution. Otherwise, we show that the energy decays to zero with rate $\frac{1}{t^{1/2}}$.

%SSSSSSSSSSSSSSSSSSSSSSSSSSSSSSSSSSSSSSSSSSSSSSSSSSSSSSSSSSSSSSSSSSSSSSSSSSSSSSSSSSSSSSSSSSSSSSSSSSSSSSSSSSSSSSSSSSSSSSSSSSSSSSSSSSSSSSSSSSSSSSSSSSSSSSSSSSSSSSSS
\section{Energy decay rate}
First, we state stability results of  system (\ref{e11})-(\ref{e16}) under one locally distributed dissipation law in the case
of equal speed wave propagation condition. We establish the following energy decay rate:

\begin{theorem} (Exponential decay rate)
Assume that $\kappa=\kappa_0$ and $\frac{\rho_1}{\rho_2}=\frac{\kappa}{b}$. Then there exist positive constants $M\geq 1$, $\omega >0$
such that for all initial data $(\varphi_0,\psi_0,\omega_0,\varphi_1,\psi_1,\omega_1)\in \HH_i,\ i=1,2,$ the energy
of the system (\ref{e11})-(\ref{e16}) satisfies the following decay
rate:
\begin{equation}
E(t)\leq Me^{-\omega t}E(0),\ \ \ \ \forall t>0.\label{ExpSta}
\end{equation}
\end{theorem}

%SSSSSSSSSSSSSSSSSSSSSSSSSSSSSSSSSSSSSSSSSSSSSSSSSSSSSSSSSSSSSSSSSSSSSSSSSSSSSSSSSSSSSSSSSSSSSSSSSSSSSSSSSSSSSSSSSSSSSSSSSSSSSSSSSSSSSSSSSSSSSSSSSSSSSSSSSSSSSSSSSSSSS

Next, under the non equal speed propagation condition, we show that the energy $E(t)$ of the
system (\ref{e11})-(\ref{e15}) loses the exponential decay rate obtained when $\kappa=\kappa_0$ and $\frac{\rho_1}{\rho_2}=\frac{\kappa}{b}$.
\begin{theorem} (Non-exponential decay rate)
Assume that $\kappa\not=\kappa_0$ or $\frac{\rho_1}{\rho_2}\not=\frac{\kappa}{b}$. Then the system (\ref{e11})-(\ref{e13}), with the boundary conditions (\ref{e15}), is not uniformly stable.
\end{theorem}
Nevertheless we establish the following polynomial type decay rate:
\begin{theorem}(Polynomial decay rate)
Assume that $\kappa\neq \kappa_0$ or $\frac{\rho_1}{\rho_2}\not=\frac{\kappa}{b}$. Then there exists a positive constant $C>0$ such that for all initial data $U_0\in D(\mathcal{A}_i)$, $i=1,2,$ the
energy of the system (\ref{e11})-(\ref{e16}) satisfies the following
decay rate:
\begin{equation}
E(t)\leq C \dfrac{1}{t} \norm{U_0}^2_{D(\mathcal{A}_j)} \ \  \ {\rm if } \ \ \kappa = \kappa_0 \  \ \ \ {\rm et } \ \ \ \ E(t)\leq C \dfrac{1}{t^{1/2}} \norm{U_0}^2_{D(\mathcal{A}_j)} \  \ \  {\rm if } \ \kappa \not= \kappa_0, \quad \forall t>0. \label{PolySta}
\end{equation}
\end{theorem}

%%%%%%%%%%%%%%%%%%%%%%%%%%%%%%%SSSSSSSSSSSSSSSSSSSSSSSSSSSSSSSSSSSSSSSSSSSSSSSSSSSSSSSSSSSSSSSSSSSSSSSSSSSSSSSSSSSSSSSSSSSSSSSSSSSSSSSSSSSSSSS

\vspace{0.5cm}
%%%%%%%%%%%%%%%%%%%%%%%%%%%%%%%%%%%%%%%%%%%%%%%%%%%%%%%%%%%%%%%%%%%%%%%%%%%%%%%%%%%%%%%%%%%%%%%%%%%%%%%%%%%%%%%%%%%%%%%%%%%%%%%%%%%%%%%%%%%%%%%%%%%%%%%%%%%
%%%%%%%%%%%%%%%%%%%%%%%%%%%%%%%%%%%%%%%%%%%%%%%%%%%%%%%%%%%%%%%%%%%%%%%%%%%%%%%%%%%%%%%%%%%%%%%%%%%%%%%%%%%%%%%%%%%%%%%%%%%%%%%%%%%%%%%%%%%%%%%%%%%%%%%%%%%%
%FFFFFFFFFFFFFFFFFFFFFFFFFFFFFFFFFFFFFFFFFFFFFFFFFFFFFFFFFFFFFFFFFFFFFFFFFFFFFFFFFFFFFFFFFFFFFFFFFFFFFFFFFFFFFFFFFFFFFFFFFFFFFFFFFFFFFFFFFFFFFFFFFFFFFFFFFFFFFFFF

\setcounter{section}{0}
\selectlanguage{francais}
\noindent{\bf II- Version fran\c{c}aise}
\section{Introduction}
Dans cette note, nous \'etudions la stabilit\'e interne locale de syst\`eme \'elastique de Bresse sous l'action d'un seul contr\^ole agissant sur l'\'equation de rotation angulaire. Nous consid\'erons le syst\`eme suivant :
\begin{eqnarray}
\rho_1\varphi_{tt}-\kappa(\varphi_x+\psi+l\omega)_x-\kappa_0 l(\omega_x-l\varphi)&=&0\hspace{1 cm} {\rm dans}\ \ (0,L)\times (0,\infty),\label{e11f}\\
\rho_2\psi_{tt}-b\psi_{xx}+\kappa(\varphi_x+\psi+l\omega)+a(x)\psi_t&=&0\hspace{1 cm} {\rm dans}\ \ (0,L)\times (0,\infty),\label{e12f}\\
\rho_1\omega_{tt}-\kappa_0(\omega_x-l\varphi)_x+\kappa l(\varphi_x+\psi+l\omega)&=&0\hspace{1cm} {\rm dans}\ \ (0,L)\times (0,\infty),\label{e13f}
\end{eqnarray}
avec les conditions aux bords de type Dirichlet-Neumann-Neumann et Dirichlet-Dirichlet-Dirichlet suivantes :
\begin{eqnarray}
\varphi(t,x)=\psi_x(t,x)=\omega_x(t,x)&=&0 \hspace{1 cm} {\rm pour}\ \ x=0,L,\label{e15f}\\
\varphi(t,x)=\psi(t,x)=\omega(t,x)&=&0 \hspace{1 cm} {\rm pour}\ \  x=0,L,\label{e16f}
\end{eqnarray}
o\`u les fonctions $\varphi$, $\psi$, $\omega$ d\'esignent, respectivement, le d\'eplacement transversal, l'angle de rotation d'un filament et le d\'eplacement
longitudinal de la poutre; la fonction $a\in L^{\infty}([0,L];\rr^+)$ v\'erifiant $a \geq a_0 >0$ dans $]\alpha,\beta[$ o\`u $0\leq \alpha< \beta\leq L$. Les coefficients $\rho_1=\rho A$, $\rho_2=\rho I$, $\kappa_0=EA$, $\kappa=\kappa^\prime GA$, $b=EI$ et $l=R^{-1}$
sont des constantes positives caract\'erisant les propri\'et\'es \'elastiques des mat\'eriaux. Physiquement, $\rho$ est la densit\'e, $E$ est le module d'\'elasticit\'e, $G$  est le module de cisaillement, $\kappa^\prime$ est le facteur de cisaillement, $A$ est l'aire de la section transversale,
$I$ est le second moment de la section transversale,
et $R$ est le rayon de courbure  (pour plus de d\'etails voir Lagnese et al. \cite{Lagnese1}).
%%%%%%%%%%%%%%%%%%%%%%%%%%%%%%%%%%%%%%%%%%%%%%%%%%%%%%%%%%%%%%%%%%%%%%%%%%%%%%%%%%%%%%%%%%%%%%%%%%%%%%%%%%%%%%%%%%%%%%%%%%%%%%%%%%%%%%%%%%%%%%%%%%%%%%%%%%%%%%%%%
%%%%%%%%%%%%%%%%%%%%%%%%%%%%%%%%%%%%%%%%%%%%%%%%%%%%%%%%%%%%%%%%%%%%%%%%%%%%%%%%%%%%%%%%%%%%%%%%%%%%%%%%%%%%%%%%%%%%%%%%%%%%%%%%%%%%%%%%%%%%%%%%%%%%%%%%%%%%%%%%%%%

Le syst\`eme de Bresse est un mod\`ele lin\'eaire couplant trois \'equations des ondes. Par ailleurs, il y a un nombre de publications r\'ecentes concernant la stabilisation de ce syst\`eme \cite{LiuRao1}, \cite{WehbeBresse}, \cite{RiveraBresse} et \cite{AlabauBresse}.
Dans \cite{AlabauBresse}, Alabau et al. ont consid\'er\'e le syst\`eme de Bresse (\ref{e11f})-(\ref{e16f}) avec une seule loi de dissipation globale agissant sur l'\'equation de rotation angulaire i.e. la fonction $a(x)$ est constante sur $]0,L[$. Dans le cas d'\'egalit\'e des vitesses de propagation ($\kappa =\kappa_0$ et $\frac{\rho_1}{\rho_2}=\frac{\kappa}{b}$), ils ont \'etabli un taux de d\'ecroissance exponentiel de l'\'energie du syst\`eme. En revanche, quand les vitesses de propagation sont diff\'erentes ($\kappa\not=\kappa_0$ ou $\frac{\rho_1}{\rho_2}\not=\frac{\kappa}{b}$), ils ont montr\'e, d'abord, que le syst\`eme n'est pas uniform\'ement stable. Ensuite, sous la condition
$\kappa =\kappa_0$ et $\frac{\rho_1}{\rho_2}\not=\frac{\kappa}{b}$, ils ont \'etabli un taux de d\'ecroissance polynomial de type $\frac{1}{t^{2/3}}$.
Puis, dans les autres cas, ils ont montr\'e que l'\'energie du syst\`eme d\'ecroit vers z\'ero comme $\frac{1}{t^{1/3}}$.
Alors, la stabilisation du syst\`eme de Bresse avec un contr\^ole interne localement distribu\'e et l'optimalit\'e du taux de d\'ecroissance polynomial restent un probl\`eme ouvert.
Notre travail consiste donc \`a \'etendre les r\'esultats de \cite{AlabauBresse}: pour prendre en consid\'eration un cas important o\`u $a$ est
une fonction dans $L^\infty([0,L];\rr_+)$ strictement positive seulement dans un intervalle ouvert quelconque $]\alpha,\beta[\subset ]0,L[$ (les cas $\alpha=0$ ou $\beta=L$ ne sont pas exclus), et pour am\'eliorer le taux de d\'ecroissance polynomial de l'\'energie.

Dans cette note, nous consid\'erons le syst\`eme de Bresse avec une seule loi de dissipation localement distribu\'ee \`a l'int\'erieur du domaine et
agissant sur l'\'equation de rotation angulaire. Sous la condition d'\'egalit\'e des vitesses de propagation
($\kappa =\kappa_0$ et $\frac{\rho_1}{\rho_2}=\frac{\kappa}{b}$),
nous \'etablissons le m\^eme taux de d\'ecroissance exponentiel de l'\'energie du syst\`eme. Dans le cas contraire,
($\kappa\not=\kappa_0$ ou $\frac{\rho_1}{\rho_2}\not=\frac{\kappa}{b}$), nous montrons, d'abord, que le syst\`eme n'est pas uniform\'ement stable. Ensuite,
dans le cas o\`u $\kappa =\kappa_0$ et $\frac{\rho_1}{\rho_2}\not=\frac{\kappa}{b}$, nous \'etablissons un nouveau taux de d\'ecroissance polynomial de type $\frac{1}{t}$. Puis, dans les autres cas, nous obtenons un nouveau taux de d\'ecroissance polynomial de type $\frac{1}{t^{1/2}}$.\\
Notons que la condition d'\'egalit\'e des vitesses de propagation n'a pas un sens physique, par cons\'equent, l'\'energie du syst\`eme d\'ecroit
vers z\'ero comme $\frac{1}{t^{1/2}}$.

%%%%%%%%%%%%%%%%%%%%%%%%%%%%%%%%%%%%%%%%%%%%%%%%%%%%%%%%%%%%%%%%%%%%%%%%%%%%%%%%%%%%%%%%%%%%%%%%%%%%%%%%%%%%%%%%%%%%%%%%%%%%%%%%%%%%%%%%%%%%%%%%%%%%%%%%
%%%%%%%%%%%%%%%%%%%%%%%%%%%%%%%%%%%%%%%%%%%%%%%%%%%%%%%%%%%%%%%%%%%%%%%%%%%%%%%%%%%%%%%%%%%%%%%%%%%%%%%%%%%%%%%%%%%%%%%%%%%%%%%%%%%%%%%%%%%%%%%%%%%%%%%%%

Soit $(\varphi,\psi,\omega)$ une solution r\'eguli\`ere du syst\`eme
(\ref{e11f})-(\ref{e13f}), avec les conditions aux bords (\ref{e15f}) ou (\ref{e16f}). Son \'energie associ\'ee est d\'efinie par :
\begin{equation}
E(t)=\frac{1}{2}\int_0^L\{\kappa\abs{\varphi_x+\psi+l\omega}^2+b\abs{\psi_x}^2+\kappa_0\abs{\omega_x-l\varphi}^2
+\rho_1|\varphi_t|^2+\rho_2|\psi_t|^2+\rho_1|\omega_t|^2\}dx. \label{e17f}
\end{equation}
Par un calcul direct nous obtenons :
\begin{equation}
\frac{d}{dt}E(t)=-\int_0^L a(x)|\psi_t|^2dx\leq 0.\label{e18f}
\end{equation}
Supposons que $L\not=\frac{n\pi}{l}$ et
posons $\Omega=(0,L)$, $L_\star^2=\{f\in L^2(0,L):\int_0^Lf(x)dx=0\}$ et $H_\star^1(\Omega)=H^1(\Omega)\cap
L_\star^2(\Omega)$. D\'efinissons les espaces d'\'energie par
$\HH_1=H_0^1\times (H_*^1)^2\times L^2\times (L_*^2)^2$ et $ \HH_2=(H_0^1)^3\times (L^2)^3$ munis de la norme :
\begin{equation}
\|U\|_{\HH_i}^2=\kappa\|\varphi_x+\psi+l\omega\|^2+b\|\psi_x\|^2+\kappa_0\|\omega_x-l\varphi\|^2
+\rho_1 \|u\|^2+\rho_2\|v\|^2+\rho_1\|z\|^2,\label{N}
\end{equation}
o\`u $\norm{\cdot}$ d\'esigne la norme dans $L^2(\Omega)$. Il est
clair que la norme (\ref{N}) est \'equivalente \`a la norme usuelle dans $\HH_i$, $i=1,2$. D\'efinissons aussi les op\'erateurs lin\'eaires $\AA_i$,
$i=1,2,$ correspondant aux conditions aux bords (\ref{e15f}) et (\ref{e16f}) respectivement, par :
\begin{equation}
D(\AA_1)=\{U\in \HH_1: \varphi\in H_0^1\cap H^2,\,\ \psi, \omega\in H_*^1\cap H^2,\,\ u, \psi_x, \omega_x\in H_0^1,\,\ v,z\in H_*^1\},\label{ope1}
\end{equation}
\begin{equation}
D(\AA_2)=\{U\in \HH_2: \varphi, \psi, \omega \in H_0^1\cap H^2,\,\ u, v, z\in H_0^1\}, \label{ope2}
\end{equation}
\begin{equation}
\AA_i(\varphi, \psi, \omega, u, v,z)=\left(\begin{array}{c}
         u\\
        v\\
        z\\
        \frac{\kappa}{\rho_1}(\varphi_x+\psi+l\omega )_x+\frac{\kappa_0 l}{\rho_1}(\omega_x-l\varphi)\\
        \frac{b}{\rho_2}\psi_{xx}-\frac{\kappa}{\rho_2 }(\varphi_x+\psi+l\omega)-\frac{1}{\rho_2}a(x)v\\
        \frac{\kappa_0}{\rho_1}(\omega_x-l\varphi)_x- \frac{\kappa l}{\rho_1}(\varphi_x+\psi+l\omega)
    \end{array}\right). \label{ope}
\end{equation}
Notons que $\AA_i$ est m-dissipatif et engendre un $C^0$-semi-groupe de contractions $e^{t\AA_i}$ sur l'espace d'\'energie $\HH_i$, $i=1,2$. Comme le
syst\`eme (\ref{e11f})-(\ref{e16f}) est \'equivalent \`a
$$
U_t= \mathcal{A}_i U \hbox{ dans } \HH_i,\ \  t>0, \quad U(0)=U_0, \ \ i=1,2
$$
avec  $U=(\varphi, \psi, \omega,\varphi_t, \psi_t, \omega_t)$, nous en  d\'eduisons l'existence et l'unicit\'e de la solution.

Comme soulign\'e pr\'ec\'edemment, nous allons \'etudier le taux de
d\'ecroissance de l'\'energie dans deux cas diff\'erents.
%SSSSSSSSSSSSSSSSSSSSSSSSSSSSSSSSSSSSSSSSSSSSSSSSSSSSSSSSSSSSSSSSSSSSSSSSSSSSSSSSSSSSSSSSSSSSSSSSSSSSSSSSSSSSSSSSSSSSSSSSSSSSSSSSSSSSSSSSSSSSSSSSSSSSSSSSSSSSSSSSSSSSSSS
\section{Stabilit\'e exponentielle et polynomiale du syst\`eme}
Dans cette section nous montrons les r\'esultats suivants :

%TTTTTTTTTTTTTTTTTTTTTTTTTTTTTTTTTTTTTTTTTTTTTTTTTTTTTTTTTTTTTTTTTTTTTTTTTTTTTTTTTTTTTTTTTTTTTTTTTTTTTTTTTTTTTT

%TTTTTTTTTTTTTTTTTTTTTTTTTTTTTTTTTTTTTTTTTTTTTTTTTTTTTTTTTTTTTTTTTTTTTTTTTTTTTTTTTTTTTTTTTTTTTTTTTTTTTTTTTTTTTTTTTTTTTTTTTT
\begin{theoreme}\label{texpdecay} (Taux de d\'ecroissance exponentiel)
Supposons que $\kappa=\kappa_0$ et $\frac{\rho_1}{\rho_2}=\frac{\kappa}{b}$. Alors il existe des constantes $M\geq 1$ et $\omega >0$ telles que pour toute
donn\'ee initiale $(\varphi_0, \psi_0, \omega_0, \varphi_1, \psi_1, \omega_1)\in \HH_j$, $j=1,2$, l'\'energie du syst\`eme (\ref{e11f})-(\ref{e16f}) satisfait l'estimation suivante :
\begin{equation}
E(t)\leq Me^{-\omega t}E(0),\ \ \ \ \forall t>0.\label{e31f}
\end{equation}
\end{theoreme}
{\it Id\'ee de la d\'emonstration}\\
D'apr\`es un r\'esultat de Huang \cite{H} et Pr\"uss \cite{Pr}, un
$C^0$-semi-groupe de contractions $e^{t\AA_j}$ sur un espace de Hilbert
$\HH_j$, $j=1,2,$ est exponentiellement stable si et seulement si les conditions suivantes :
$${\rm (H1)}:\ \ \  i\rr\,\subset\,\rho(\mathcal{A}_j) \ \ \ \ \ {\rm et } \ \ \ \ \  {\rm (H2)}:\ \ \ \sup_{\lambda \in\rr}\norm{(i\lambda
I-\mathcal{A}_j)^{-1}}<\,+\infty$$
sont v\'erifi\'ees. Nous montrons que $\AA_j$ n'a
pas de valeur propre sur l'axe imaginaire. Alors, le fait que la r\'esolvante de $\AA_j$ est compacte,
entraine bien la condition (H1). La condition (H2) se d\'emontre en utilisant un argument de contradiction (voir th\'eor\`eme(\ref{theo})). \hfill$\Box$
%TTTTTTTTTTTTTTTTTTTTTTTTTTTTTTTTTTTTTTTTTTTTTTTTTTTTTTTTTTTTTTTTTTTTTTTTTTTTTTTTTTTTTTTTTTTTTTTTTTTTTTTTTTTTTTTTTTTTTTTT
\begin{theoreme}(Taux non-exponentiel)
Supposons que $\kappa\neq \kappa_0$ ou $\frac{\rho_1}{\rho_2}\not=\frac{\kappa}{b}$. Alors, le syst\`eme (\ref{e11f})-(\ref{e13f}) avec les conditions aux bords (\ref{e15f}), n'est pas uniform\'ement stable.\\
\end{theoreme}
{\it Id\'ee de la d\'emonstration}\\
Nous construisons une suite $(\lambda_n)\subset \rr$ et une suite $(U^n)\subset D(\AA_1)$, telles que
$\abs{\lambda_n}\longrightarrow+\infty$, $\norm{U^n}_{{\HH}_1}=1$ et
$\norm{(i\lambda_nI-\mathcal{A}_1)U^n}_{{\HH}_1}\longrightarrow 0$. Alors, la r\'esolvante de
$\AA_1$ n'est pas uniform\'ement born\'ee sur l'axe imaginaire et par cons\'equent le syst\`eme (\ref{e11f})-(\ref{e15f}) n'est pas uniform\'ement stable
(voir \cite{H} et \cite{Pr}). \hfill$\Box$
%SSSSSSSSSSSSSSSSSSSSSSSSSSSSSSSSSSSSSSSSSSSSSSSSSSSSSSSSSSSSSSSSSSSSSSSSSSSSSSSSSSSSSSSSSSSSSSSSSSSSSSSSSSSSSSSSSSSSSSSSSSSSSSSSSSSSSSSSSSSSSSSSSSSSS
%TTTTTTTTTTTTTTTTTTTTTTTTTTTTTTTTTTTTTTTTTTTTTTTTTTTTTTTTTTTTTTTTTTTTTTTTTTTTTTTTTTTTTTTTTTTTTTTTTTTTTTTTTTTTTTTTTTTTTTTTTTTTTTTTTTTTTTTTTTTTTTTTTTTT
%%%%%%%%%%%%%%%%%%%%%%%%%%%%%%%%%%%%%%%%%%%%%%%%%%%%%%%%%%%%%%%%%%%%%%%%%%%%%%%%%%%%%%%%%%%%%%%%%%%%%%%%%%%%%%%%%%%%%%%%%%%
\begin{theoreme}(Taux polynomial)\label{theo}
Supposons que $\kappa\neq \kappa_0$ ou $\frac{\rho_1}{\rho_2}\not=\frac{\kappa}{b}$. Alors il existe une constante $C>0$ telle que pour toute donn\'ee initiale
$U_0=(\varphi_0, \psi_0, \omega_0, \varphi_1, \psi_1, \omega_1)\in D(\mathcal{A}_j)$, j=1,2, l'\'energie du syst\`eme
(\ref{e11f})-(\ref{e13f}), avec (\ref{e15f}) ou (\ref{e16f}), satisfait l'estimation suivante :
\begin{equation}
E(t)\leq C \dfrac{1}{t} \norm{U_0}^2_{D(\mathcal{A}_j)} \ \  \ {\rm si } \ \ \kappa = \kappa_0 \  \ \ \ {\rm et } \ \ \ \ E(t)\leq C \dfrac{1}{t^{1/2}} \norm{U_0}^2_{D(\mathcal{A}_j)} \  \ \  {\rm si } \ \kappa \not= \kappa_0, \quad \forall t>0. \label{e41}
\end{equation}
\end{theoreme}
%%%%%%%%%%%%%%%%%%%%%%%%%%%%
On suppose que $\kappa \not= \kappa_0$.
La d\'emonstration du th\'eor\`eme (\ref{theo}) n\'ecessite la d\'emonstration de plusieurs lemmes.\\ \\
Soit  $(\lambda_n)$ une suite dans $\rr$ et $(U^n)$ une suite dans $D(\AA_j)$, telles que
\begin{equation}
\abs{\lambda_n}\longrightarrow+\infty, \ \ \ \norm{U^n}_{{\HH}_j}=\norm{(\varphi^n, \psi^n, \omega^n,u^n,v^n,z^n)}_{{\HH}_j}=1,\label{Ubornee}
\end{equation} et
\begin{equation}
\lambda_n^{4}(i\lambda_nI-\mathcal{A}_j)(\varphi^n, \psi^n, \omega^n,u^n,v^n,z^n)=(f^n_1,f_2^n,f_3^n,g^n_1,g^n_2,g^n_3)\longrightarrow
0 \quad{\rm dans }\quad\mathcal{H}_j.\label{e33}
\end{equation}
Nous \'ecrivons (\ref{e33}) sous la forme :
\begin{eqnarray}
\lambda_n^2\varphi^n+\frac{\kappa}{\rho_1}(\varphi_{xx}^n+\psi_x^n+l\omega^n_x)+\frac{\kappa_0 l}{\rho_1}(\omega_x^n-l\varphi^n)&=&-\frac{g^n_1+i\lambda_nf^n_1}{\lambda_n^{4}},
\label{fi}\\
\lambda_n^2\psi^n+\dfrac{b}{\rho_2}\psi^n_{xx}-\dfrac{\kappa}{\rho_2}(\varphi^n_x+\psi^n+l\omega^n)-\dfrac{1}{\rho_2 }a(x)v^n&=&-\frac{g^n_2+i\lambda_nf^n_2}{\lambda_n^{4}}
,\label{psi}\\
\lambda_n^2\omega^n+\dfrac{\kappa_0}{\rho_1}(\omega^n_{xx}-l\varphi^n_x)-\dfrac{\kappa l}{\rho_1}(\varphi^n_x+\psi^n+l\omega^n)
&=&-\frac{g^n_3+i\lambda_nf^n_3}{\lambda_n^{4}}.\label{omega}
\end{eqnarray}
%%%%%%%%%%%%%%%%%%%%%%%%%%%%%%%%%%%%%%%%%%%%%%%%%%%%%%%%%%%
%EEEEEEEEEEEEEEEEEEEEEEEEEEEEEEEEEEEEEEEEEEEEEEEEEEEEEEEEEEEEEEEEEEEEEEEEEEEEEEEEEEEEEEEEEEEEEEEEEEEEEEEEEEEEEEEEEEEEEEEEEEEEEEEEEEEEEEEEEEEEEEEEEEEEEE
\begin{lemme}\label{lem1}
Sous les m\^emes notations pr\'ec\'edentes, nous avons
\begin{equation}
\int_0^L a(x)\abs{\psi^n}^2dx=\frac{o(1)}{\lambda_n^{6}}.
\end{equation}
\end{lemme}
{\it D\'emonstration}\\
Il suffit de multiplier (\ref{e33}) par $U^n=(\varphi^n, \psi^n, \omega^n,u^n,v^n,z^n)$ pour obtenir :
\begin{equation}
\int_0^L a(x) \abs{v^n}^2dx={\rm Re} ((i\lambda_n I-\mathcal{A}_j)U^n,U^n)_{\mathcal{H}_j}=\frac{o(1)}{\lambda_n^{4}}\ \ \ \ {\rm et} \  \
\int_0^L a(x)\abs{\psi^n}^2dx=\frac{o(1)}{\lambda_n^{6}}.\label{dissipation} \hfill\Box
\end{equation}
%EEEEEEEEEEEEEEEEEEEEEEEEEEEEEEEEEEEEEEEEEEEEEEEEEEEEEEEEEEEEEEEEEEEEEEEEEEEEEEEEEEEEEEEEEEEEEEEEEEEEEEEEEEEEEEEEEEEEEEEEEEEEEEEEEEEEEEEEEEEEEEEEEEEEEE
%%%%%%%%%%%%%%%%%%%%%%%%%%%%%%%%%%%%%%%%%%%%%%%%%%%%%%%%%%%%%%%%%%%%%%%%%%%%%%%%%%%%%%%%%%%%%%%%%%%%%%%%%%%%%%%%%%%%%%%%%%%%%%%%%%%%%%%%%%%%%%%%%%%%%
Soit $\varepsilon >0$ tel que $\alpha<\alpha + \varepsilon <\beta-\varepsilon <\beta$.
Consid\'erons une fonction $\eta \in C^1([0,L])$ d\'efinie par $0\leq\eta\leq 1$, $\eta=1$ sur $[\alpha+\varepsilon, \beta-\varepsilon]$ et $\eta=0$ sur
$[0, \alpha] \cup [\beta,L]$ .

%EEEEEEEEEEEEEEEEEEEEEEEEEEEEEEEEEEEEEEEEEEEEEEEEEEEEEEEEEEEEEEEEEEEEEEEEEEEEEEEEEEEEEEEEEEEEEEEEEEEEEEEEEEEEEEEEEEEEEEEEEEEEEEEEEEEEEEEEEEEEEEEEEEEEEE
\begin{lemme}\label{lem2}
Sous les m\^emes notations pr\'ec\'edentes, nous avons
\begin{equation}
\int_0^L \eta\abs{\psi_x^n}^2=\frac{o(1)}{\lambda_n^{3}}.\label{infopsix1}
\end{equation}
\end{lemme}
{\it D\'emonstration}\\
Dans un premier temps, notons que d'apr\`es (\ref{Ubornee}) et (\ref{e33}) nous avons $\norm{\varphi^n}=O(1/\lambda_n)$, $\norm{\psi^n}=O(1/\lambda_n)$ et $\norm{\omega^n}=O(1/\lambda_n)$.
Ensuite, multiplions l'\'equation (\ref{psi}) par $\eta \bar{\psi}^n$, nous obtenons :
\begin{equation}
b\int_0^L \eta \abs{\psi^n_x}^2=\rho_2\int_0^L \eta \abs{\lambda_n\psi^n}^2-b\int_0^L \eta^\prime \psi^n_x\bar{\psi^n}-
\int_0^L [\kappa(\varphi_x^n+\psi^n +l\omega^n)+a v^n]\eta\bar{\psi}^n
+\frac{o(1)}{\lambda_n^{4}}.\label{psipsix}
\end{equation}
Nous utilisons (\ref{dissipation})  et le fait que
$\norm{\psi^n_x}$, $\norm{\varphi_x^n+\psi^n +l\omega^n}$ sont uniform\'ement born\'ees dans (\ref{psipsix})$ \times \lambda_n^3$. Nous d\'eduisons alors (\ref{infopsix1}). \hfill$\Box$
%EEEEEEEEEEEEEEEEEEEEEEEEEEEEEEEEEEEEEEEEEEEEEEEEEEEEEEEEEEEEEEEEEEEEEEEEEEEEEEEEEEEEEEEEEEEEEEEEEEEEEEEEEEEEEEEEEEEEEEEEEEEEEEEEEEEEEEEEEEEEEEEEEEEEEEEEEEEEEEEEEEEEEE
%EEEEEEEEEEEEEEEEEEEEEEEEEEEEEEEEEEEEEEEEEEEEEEEEEEEEEEEEEEEEEEEEEEEEEEEEEEEEEEEEEEEEEEEEEEEEEEEEEEEEEEEEEEEEEEEEEEEEEEEEEEEEEE

%LEMMMMMMMMMMMMMMMMEEEEEEEEEEEEEEEEEEEEEEEEEEEEEEEEEEEEEEEEEEEEEEEEEEEEEEEEEEEEEEEEEEEEEEEEEEEEEEEEEEEEEEEEEEEEEEEEEEEEE
\begin{lemme}
Sous les m\^emes notations pr\'ec\'edentes, nous avons
\begin{equation}
\int_0^L\eta|\varphi_x^n|^2=o(1)\hspace{1cm} {\it et} \hspace{1cm}\int_0^L\eta|\varphi^n|^2=\frac{o(1)}{\lambda_n^2}.\label{fix}
\end{equation}
\end{lemme}
{\it D\'emonstration}\\
Multiplions (\ref{psi}) par $\eta\bar{\varphi}^n_x$, (\ref{fi}) par $\eta\bar{\varphi}^n$ et utilisons (\ref{dissipation}) et le fait que $\norm{\varphi_x^n}=O(1)$. Nous obtenons respectivement :
\begin{equation}
\kappa\int_0^L \eta \abs{\varphi^n_x}^2=\rho_2\int_0^L \eta\lambda_n^2 \psi^n \bar{\varphi_x^n} -b \int_0^L \eta\lambda_n\psi^n_x\lambda_n^{-1}\bar{\varphi^n}_{xx}
-b \int_0^L\eta^\prime\psi^n_x\bar{\varphi^n}_{x}
- \int_0^L (\kappa l\omega^n+a v)\eta\bar{\varphi}^n_x
%-\kappa\int_0^L \eta\psi^n\bar{\varphi}^n_x -\int_0^L av^n\eta\bar{\varphi}^n_x
+\frac{o(1)}{\lambda_n^{3}}\label{fi2}
\end{equation}
et
\begin{equation}
\rho_1\int_0^L \eta \abs{\lambda_n\varphi^n}^2=\kappa\int_0^L (\eta \abs{\varphi^n_x}^2+ (\eta^\prime \varphi_x-\eta\psi^n_x)\bar{\varphi^n})
+l\int_0^L(\kappa+\kappa_0)\omega^n(\eta\bar{\varphi^n})_x+l^2\kappa_0 \int_0^L\eta \abs{\varphi^n}^2 +\frac{o(1)}{\lambda_n^{4}}.\label{fi3}
\end{equation}
D'abord, utilisons (\ref{dissipation}), (\ref{infopsix1}) et le fait que $\norm{\varphi_x^n}=O(1)$, $\norm{\varphi_{xx}^n}=O(\lambda_n)$ et $\norm{\omega^n}=o(1)$ dans (\ref{fi2}), nous obtenons la premi\`ere partie de (\ref{fix}). Ensuite, utilisons de plus le fait que $\norm{\varphi^n}=\norm{\varphi_x^n}=o(1)$ dans (\ref{fi3}), nous d\'eduisons la deuxi\`eme. \medskip\hfill$\Box$

%%%%%%%%%%%%%%%lemme Lemme Lemme Lemmme Lemme Lemme
\begin{lemme}\label{lalemme}
Soit $\frac{1}{2} \leq\gamma\leq 1 $. Sous les m\^emes notations pr\'ec\'edentes, supposons que
\begin{equation}
\displaystyle\int_0^L \eta\abs{\psi_x^n}^2=\dfrac{o(1)}{\lambda_n^{2+2\gamma}} \label{hyplem}.
\end{equation}
Alors, nous avons
\begin{equation}
\displaystyle\int_0^L \eta\abs{\varphi_x^n}^2=\frac{o(1)}{\lambda_n^{2\gamma}} \label{reslem}.
\end{equation}
\end{lemme}
{\it D\'emonstration}\\
Soit $l_N=\displaystyle\sum_{k=0}^{N}\frac{1}{2^k}$, nous allons d\'emontrer par r\'ecurrence sur $N$ que
$\displaystyle\int_0^L \eta\abs{\varphi_x^n}^2=\dfrac{o(1)}{\lambda_n^{\gamma l_N}}.$\\
Utilisons (\ref{dissipation}), (\ref{fix}), (\ref{hyplem}) et le fait que $\norm{\varphi_{xx}^n}=O(\lambda_n)$ et $\norm{\omega^n}=O(\frac{1}{\lambda_n})$ dans (\ref{fi2})$\times\lambda_n^\gamma$, nous obtenons
$\displaystyle\int_0^L \eta|\varphi_x^n|^2=\frac{o(1)}{\lambda_n^\gamma}$. Alors, l'\'egalit\'e est vraie pour $N=0$. Supposons
\begin{equation}
\displaystyle\int_0^L\eta|\varphi_x^n|^2=\frac{o(1)}{\lambda_n^{\gamma l_N}}. \label{hyprec}
\end{equation}
Utilisons (\ref{fix}), (\ref{hyplem}), (\ref{hyprec}) et le fait que $\norm{\omega^n}=O(\frac{1}{\lambda_n})$ et $\gamma l_N\leq 2$ dans (\ref{fi3})$\times\lambda_n^{\gamma l_N}$, nous obtenons
\begin{equation}
\displaystyle\int_0^L\eta|\varphi^n|^2=\frac{o(1)}{\lambda_n^{2+\gamma l_N}}.\label{resrec}
\end{equation}
Utilisons (\ref{resrec}) dans (\ref{fi}) nous obtenons
\begin{equation}
\norm{\eta\varphi_{xx}^n}=O(\lambda_n^{1-\frac{\gamma}{2}l_N}).\label{fixx}
\end{equation}
Notons que $\gamma + \frac{\gamma}{2}l_N= \gamma l_{N+1}$.
En utilisant (\ref{dissipation}), (\ref{hyplem}), (\ref{hyprec}), (\ref{resrec}), (\ref{fixx}) et le fait que $\norm{\omega^n}=O(\frac{1}{\lambda_n})$
dans (\ref{fi2})$\times \lambda_n^{\gamma l_{N+1}}$, nous obtenons
$\displaystyle\int_0^L \eta\abs{\varphi_x^n}^2=\dfrac{o(1)}{\lambda_n^{\gamma l_{N+1}}}.$
Par cons\'equent, (\ref{hyprec}) est vraie pour tout $N\geq 0$.
Enfin, en utilisant le fait que $\displaystyle\sum_{k=0}^{+\infty}\frac{1}{2^k}=2$, nous d\'eduisons (\ref{reslem}).\hfill$\Box$
%%%%%%%%%%%%%%%%%%%%%%%%%%%%%%%%%%%%%%%%%%%%%%%%%%%%%%%%%%%%%%%%%%%%%%%%%%%%%%%%%%%%%%%%%%%%%%%%%%%%%%%%%%%%%%%%%%%%%%%%%%%%%%%%%%%%%%%%%%%%%%%%%
%%%%%%%%%%%%%%EEEEEEEEEEEEEEEEEEEEEEEEEEEEEEEEEEEEEEEEEEEEEEEEEEEEEEEEEEEEEEEEEEEEEEEEEEEEEEEEEEEEEEEEEEEEEEEEEEEEEEEEEEEEEEEEEEEEEEEEEEEEEEEEEEE
\begin{lemme}\label{lem5}
Sous les m\^emes notations pr\'ec\'edentes, nous avons
\begin{equation}
\int_0^L \eta\abs{\psi_x^n}^2=\frac{o(1)}{\lambda_n^{4}},\ \ \ \ \int_0^L \eta\abs{\varphi_x^n}^2=\frac{o(1)}{\lambda_n^{2}}\ \ \ \ {\rm et}\ \ \ \ \int_0^L \eta\abs{\varphi^n}^2=\frac{o(1)}{\lambda_n^{4}}.\label{le res}
\end{equation}
\end{lemme}
{\it D\'emonstration}\\
Soit $N\in \nn^\star$, posons $\hat{l}_N=\displaystyle\sum_{k=1}^{N}\dfrac{1}{2^k}$. Nous allons d\'emontrer par r\'ecurrence sur $N$ que
$\displaystyle \int_0^L \eta \abs{\psi_x^n}^2=\frac{o(1)}{\lambda_n^{2+2\hat{l}_N}}.$ D'abord, l'\'equation (\ref{infopsix1}) assure l'\'egalit\'e pour $N=1$.
Ensuite, supposons que
\begin{equation}
\int_0^L \eta\abs{\psi_x^n}^2=\frac{o(1)}{\lambda_n^{2+2\hat{l}_{N-1}}}.\label{psixN-1}
\end{equation}
Notons que $3+\hat{l}_{N-1}=2+2\hat{l}_N$ et utilisons (\ref{dissipation}), (\ref{psixN-1}), le lemme(\ref{lalemme}) et le fait que $\norm{\omega^n}=O(\frac{1}{\lambda_n})$ dans $(\ref{psipsix})\times\lambda_n^{2+2\hat{l}_N}$. Nous d\'eduisons alors le r\'esultat pour tout entier $N\geq 1$.
Finalement, en utilisant le fait que $\displaystyle\sum_{k=1}^{+\infty}\dfrac{1}{2^k}=1$, le lemme (\ref{lalemme}) et l'\'equation (\ref{fi3}),
nous concluons (\ref{le res}).\hfill$\Box$
%EEEEEEEEEEEEEEEEEEEEEEEEEEEEEEEEEEEEEEEEEEEEEEEEEEEEEEEEEEEEEEEEEEEEEEEEEEEEEEEEEEEEEEEEEEEEEEEEEEEEEEEEEEEEEEEEEEEEEEEEEEEEEEEEEEEE

\begin{lemme}\label{lem6}
Sous les m\^emes notations pr\'ec\'edentes, nous avons
\begin{equation}
\int_0^L \eta\abs{\omega_x^n}^2= o(1).  \label{omega3}
\end{equation}
\end{lemme}
{\it D\'emonstration}\\
Multiplions l'\'equation (\ref{fi}) par $\eta\bar{\omega^n_x}$ et utilisons (\ref{le res})  et le fait que $\norm{\omega_x^n}$ est uniform\'ement born\'ee.
Nous d\'eduisons :
\begin{equation}
(\kappa+\kappa_0)l\int_0^L\eta\abs{\omega_x^n}^2= \kappa\int_0^L \eta\lambda_n \varphi^n_x \lambda_n^{-1}\bar{\omega^n_{xx}}
+o(1).\label{omegax}
\end{equation}
Utilisons (\ref{le res}) et le fait que $\norm{\omega_{xx}^n}=O(\lambda_n)$ dans (\ref{omegax}),
nous concluons (\ref{omega3}).\hfill$\Box$\\ \\
%%%%%%%%%%%%%%%%%%%%%%%%%%%%%%%%EEEEEEEEEEEEEEEEEEEEEEEEEEEEEEEEEEEEEEEEEEEEEEEEEEEEEEEEEEEEEEEEEEEEEEEEEEEEEEEEE%%%%%%%%%%%%%%%%%%%%%%%%%%%%%%%%%%%%%%%%%%%%%%%%%%
%%%%%%%%%%%%%%%%%%%%%%%%%%%%%%%%%%%%%%%%%%%%%%%%%%%%%%%%%%%%%%%%%%%%%%%%%%%%%%%%%%%%%%%%%%%%%%%%%%%%%%%%%%%%%%%%%%%%%%%%%%%%%%%%%%%%%%%%%%%%%%%%%%%%%%%%%%%%%%%
{\bf  D\'emonstration du th\'eror\`eme(\ref{theo})} (cas $\kappa\not=\kappa_0)$\\
%%Sous les m\^emes notations pr\'ec\'edentes, nous avons
%%\begin{equation}
 %%\norm{U^n}_{{\HH}_j}=o(1) \ \ sur\ \ ]0,\alpha+\varepsilon[\ \cup\ ]\beta-\varepsilon,L[.\label{res}
%%\end{equation}
Nous d\'efinissons une fonction $\hat{\eta} \in C^1([0,L])$ v\'erifiant $0\leq\hat{\eta}\leq 1$, $\hat{\eta}=1$ sur $[0,\alpha+\varepsilon]$ et $\hat{\eta}=0$ sur
$[\alpha+2\varepsilon, L]$.\\
Multiplions (\ref{fi}) par $2\rho_1 x\hat{\eta}\bar{\varphi^n_x}$, (\ref{psi}) par $2\rho_2 x\hat{\eta}\bar{\psi^n_x}$ et (\ref{omega}) par $2\rho_1 x\hat{\eta}\bar{\omega^n_x}$ et sommons les identit\'es.
En utilisant le fait que $\norm{U^n}_{{\HH}_j}=o(1)\ \  sur\ \ ]\alpha+\varepsilon, \beta-\varepsilon[$ (lemmes (\ref{lem5}) et (\ref{lem6}))
%\begin{equation}
 %\norm{\eta\varphi_x^n}=\norm{\eta\psi_x^n}=\norm{\eta\omega_x^n}=
% \norm{\eta\lambda_n\varphi^n}=\norm{\eta\lambda_n\psi^n}=
%\norm{\eta\lambda_n\omega^n}=o(1), \label{res1}
%\end{equation}
et le fait que
${\rm supp}\hat{\eta}^\prime\subseteq {\rm supp}\eta,$
nous obtenons
%\begin{equation}
%\rho_1\int_0^L \hat{\eta} \abs{\lambda_n\varphi^n}^2+
%\rho_2\int_0^L \hat{\eta} \abs{\lambda_n\psi^n}^2+
%\rho_1\int_0^L \hat{\eta} \abs{\lambda_n\omega^n}^2 +
%\kappa\int_0^L \hat{\eta} \abs{\varphi_x^n}^2+
%b\int_0^L \hat{\eta} \abs{\psi_x^n}^2+\kappa_0\int_0^L \hat{\eta} \abs{\omega_x^n}^2=o(1).
%\end{equation}
%D'o\`u
$\norm{U^n}_{{\HH}_j}=o(1)\ \  sur\ \ ]0,\alpha+\varepsilon[.$
Par le m\^eme raisonnement sur $]\beta-\varepsilon,L[$, nous concluons que $\norm{U^n}_{{\HH}_j}=o(1) \ \ sur\ \ ]\beta-\varepsilon,L[$.
Par cons\'equent, nous avons
$\norm{U^n}_{{\HH}_j}=o(1)$, ce qui contredit le fait que $\norm{U^n}_{{\HH}_j}=1$ et nous permet de d\'eduire que
$$\sup_{\lambda \in\rr}\frac{1}{\abs{\lambda}^{4}}\norm{(i\lambda I-\mathcal{A}_j)^{-1}}<\,+\infty.$$
Finalement, sous le fait que $i\rr\,\subset\,\rho(\mathcal{A}_j)$ (th\'eor\`eme (\ref{texpdecay})), le "th\'eor\`eme 2.4" de
\cite{BT} (voir aussi  \cite{LiuRao2,Batty}) affirme que l'\'energie du syst\`eme d\'ecroit comme $\frac{1}{t^{\frac{1}{2}}}$.
\hfill$\Box$\medskip

\bibliographystyle{abbrv} %Biblionote
\bibliography{biblio-CrasBresse} %Biblionote

\end{document}